\documentclass[reqno]{amsart}
\usepackage[applemac]{inputenc}

\usepackage{amsmath}
\usepackage{amssymb}
\usepackage{amsfonts}
\usepackage{amsthm}
\usepackage{multirow}
\usepackage{enumerate}
\usepackage{color}
\usepackage[normalem]{ulem}%\xout{text to be crossed out}\sout{text to be striked out} \usepackage{cancel}

\usepackage{array}
\newcolumntype{L}{>{\displaystyle}l}
\newcolumntype{C}{>{\displaystyle}c}
\newcolumntype{R}{>{\displaystyle}r}

\usepackage{setspace}
\newcommand{\R}{\mathbb{R}}
\newcommand{\C}{\mathbb{C}}
\newcommand{\f}{\rightarrow}
\newcommand{\deb}{\bar\partial}
\newcommand{\de}{\partial}
\newcommand{\K}{K\"{a}hler}

\newcommand{\lmb}{\lambda}
\newcommand{\ov}[1]{\overline{#1}}
\newcommand{\w}[1]{\widetilde{#1}}

\newcommand{\D}{\mathcal{D}}

\newcommand{\hyp}{h}

\newcommand{\tr}{\operatorname{tr}}
\newcommand{\Tr}{\operatorname{Tr}}

\newcommand{\arctanh}{\operatorname{arctanh}}

\newcommand{\Aut}{\operatorname{Aut}}
\newcommand{\U}{\operatorname{U}}
\newcommand{\Vol}{\operatorname{Vol}}

\newcommand{\B}{\operatorname{\mathcal B}}
\newcommand{\rank}{\operatorname{rank}}

\newcommand{\grad}{\operatorname{grad}}
\newcommand{\Ent}{\operatorname{Ent_d}}

\newcommand{\Jac}{\operatorname{Jac}}
\newcommand{\Bary}{\operatorname{Bar}}
\newcommand{\Isom}{\operatorname{Isom}}

\newcommand{\ep}{{\varepsilon}}
\newcommand{\X}{{\mathcal X}}

\newcommand{\Entv}{\operatorname{Ent_v}}
\newcommand{\di}{\rho}
\newcommand{\Hess}{\operatorname{Hess}}

\newtheorem{thm}{Theorem}%[section]
\newtheorem{prop}{Proposition}[section]
\newtheorem{lem}[prop]{Lemma}
\newtheorem{cor}[prop]{Corollary}
\newtheorem{defn}[prop]{Definition}

\newtheorem{rmk}[prop]{Remark}

\newtheorem*{theorA}{Theorem A}
\newtheorem*{theorB}{Theorem B}
\newtheorem*{theorC}{Theorem C}

\begin{document}

\date{\today}

\author[R. Mossa]{Roberto Mossa}

\address{Departamento de Matem\'atica
\newline\indent
Universidade de S\~ao Paulo
\newline\indent
Rua do Mata\~o 1010, 05508-090 S\~ao Paulo, SP, Brazil}
\email{robertom@ime.usp.br}

\thanks{The author research was supported by
FAPESP (FAPESP grant: 2018/08971-9)}

%\title[A gap theorem for   maps into complex hyperbolic manifolds]
%{A gap theorem for  continuous maps into complex hyperbolic manifolds through diastatic entropy}

\title[On the diastatic entropy and $\mathcal C^1 $-rigidity of complex hyperbolic manifolds]
{On the diastatic entropy and $\mathcal C^1 $-rigidity of complex hyperbolic manifolds}

\begin{abstract}
Let $f: (Y, g) \f (X,g_0)$  be a non zero degree continuous map between compact \K\ manifolds of dimension $n\geq 2$, where $g_0$ has constant negative holomorphic sectional curvature. Adapting the Besson--Courtois--Gallot \emph{barycentre map} techniques to the \K\ setting, we prove a   gap theorem  in terms of the degree of $f$ and the diastatic entropies of $(Y, g)$ and $(X,g_0)$ which extends the rigidity result proved by the author in \cite{M06}.

\end{abstract}

\maketitle

\tableofcontents

\section{Introduction and statement of the main results}
It is a classical problem to determine when a continuous map between two closed smooth manifolds is homotopic to a more regular one. Of course, the father of this problems is the celebrated Mostow Rigidity Theorem which was beautifully extended in the seminal paper \cite{bcg1} (see also \cite{bcg0, bcg2, bcg3}) by G. Besson, G. Courtois and S. Gallot. This is  expressed by the following result which combined with  barycentre techniques developed in its proof has provided a solution of long-standing problems. {Denoted by $\Entv(M,g)$ the volume entropy of a compact Riemannian manifold $(M,g)$ we have:} 
\begin{theorA}[G. Besson, G. Courtois, S. Gallot]\label{thm dentropy}
Let  $(Y, g)$ be a compact Riemannian manifold of dimension $n\geq 3$ and let $(X,g_0)$ be a compact negatively curved locally symmetric Riemannian manifold of the same dimension of $Y$.
If $f:Y \rightarrow X$ is a nonzero degree continuous map, then
\begin{equation}\label{eq main}%\small
\Entv^{n}\left(Y,\, g \right)\, {\Vol}\left(Y,\, g \right)\,\geq\, \left|{\deg} \left( f \right)\right| \, \Entv^{n}\left(X,\, g_0\right)\, {\Vol}\left(X,\, g_0\right).
\end{equation}
Moreover, the equality is attained if and only if $f$ is homotopic to a {homothetic}  covering $F:Y \f X$.
\end{theorA}
The following theorem (Theorem B), proved by the author of the present paper in  \cite{M06}, represents an  extension of Theorem A in the \K\ setting
by substituting  the volume entropy with the diastatic entropy (introduced in \cite{mossau} and studied in  \cite{mossal} in the homogeneous setting).
We briefly recall its definition for reader convenience.
Let  $\pi : \left( \w Y,\, \w g \right) \f  \left(Y,\, g \right)$ be the universal \K\ covering (i.e. $\pi$ is a holomorphic covering map and $\pi^* g = \w g$) of a compact \K\ manifold $\left( Y,\, g \right)$ and assume that the diastasis function $\D: \w Y \times \w Y \f \R$ is \emph{globally defined}, that is, defined in whole $\w Y \times \w Y$ (see next section of the definition of diastasis function). Then, the {diastatic entropy} of $\left( Y, \, g \right)$ is the \K\ invariant of $g$ given by
\begin{equation}\label{int entrop}
\begin{split}
\Ent\left(Y,\,  g\right)=\operatorname{ \mathcal X} \left({\w g}\right)\,\inf\left\{c \in \R^+: \,  \int_{\w Y} e^{-c\, \D_{q}}\, \nu_{\w g} < \infty\right\},
\end{split}
\end{equation}
where $\mathcal X \left({\w g} \right)= \sup_{p,\, q\,\in\, \w Y}{\|\grad_p \D_q\|}$ and $\nu_{\w g}$ is the volume form associated to ${\w g}$. If $\mathcal X \left({\w g}\right)=\infty$ or the infimum in \eqref{int entrop} is not achieved by any $c \in \R^+$, we set $\Ent\left(Y,\, g\right)=\infty$.
It is not hard to see that this definition is independent on the  point $q$ (see \cite{M06} for details).

\begin{theorB}\label{thm dentropy}
Let $\left(Y,\, g\right)$ be a compact \K\ manifold of complex dimension $n\geq 2$ and let $\left(X,\, g_0\right)$ be a compact complex hyperbolic manifold\footnote{Notice that  a negatively curved locally  hermitian symmetric \K\ manifold is authomatically a complex hyperbolic manifold, namely its holomorphic sectional curvature is constant. This is the reason, together with the use of diastatic entropy instead of the  volume entroopy,  why Theorem B can be considered an extension to the \K\ setting of Theorem A.}  of the same dimension of $Y$. If $f:Y \rightarrow X$ is a nonzero degree continuous map, then
\begin{equation}\label{eq main}%\small
\Ent^{2n}\left(Y,\, g \right)\, {\Vol}\left(Y,\, g \right)\,\geq\, \left|{\deg} \left( f \right)\right| \, \Ent^{2n}\left(X,\, g_0\right)\, {\Vol}\left(X,\, g_0\right).
\end{equation}
Moreover, if $g$ and $g_0$ are rescaled so that
$
\Ent\left( Y,\, g\right) =\Ent \left( X,\, g_{0}\right),
$
 the equality is attained if and only if $f$ is homotopic to a holomorphic or anti-holomorphic   isometric covering $F: Y \f X$.
\end{theorB}
Later on, S. Gallot extends Theorem A by proving the  following gap result (Theorem C).
Before stating his result we need  the following definitions.
We say that  a Riemannian manifold  $\left(Y,\, g\right)$ of dimension $m$ has {\em bounded Hessian} if, for any point  $p$
of  its  Riemannian universal covering $\left(\w Y,\, \w g\right)$, there exists a positive constant $C$ %{\color{red} independent of p,}
such that   $\left| \lmb^p_j\right| < C$, for all  $j=1,\dots,m$, where $\lmb^p_j$ are the eingenvalues
of  the  Hessian  of $\tilde\rho_p$, the geodesic distance {from $p$}.
We say that  a family   $F_\epsilon :Y\rightarrow X$, $\epsilon >0$,  of $\mathcal {C}^1$-{maps}   between two compact Riemannian manifolds of the same dimension $m$  is  {\em almost-isometric}
if there exist two constants $A'\left(\epsilon\right)$  and $A''\left(\epsilon\right)$ determined by $m$ and $\epsilon$ such that
\begin{equation*}
A'\left(\epsilon\right) \leq \frac{\left\|d_p F_\epsilon\left(u\right)\right\|}{\left\|u\right\|}\leq A''\left(\epsilon\right) 
\end{equation*}
where $A'\left(\epsilon\right) \f 1$ and $A''\left(\epsilon\right) \f 1$ as $\epsilon \f 0$.

\begin{theorC}[S. Gallot (unpublished, private comunications)]\label{thm gallot}Let  $(Y, g)$ be a compact Riemannian manifold {with bounded Hessian} of dimension $m\geq 3$ and let $(X,g_0)$ be a compact negatively curved locally symmetric Riemannian manifold  of the same dimension of $Y$.
If $f:Y \f X$ is a non zero degree continuous map  and there exists   a sufficiently small positive constant $\epsilon$ such that
\begin{equation*}\label{EQ perto vol ent}
\Entv^m\left( Y,\, g\right)\, \Vol\left( Y,\, g\right) - \left| \deg \left( f \right) \right| \, \Entv ^m\left( X,\, g_{0}\right)\, \Vol \left( X,\, g_{0}\right)  \leq \epsilon,
\end{equation*}
then $f$ is homotopic to a $\mathcal {C}^1$-{covering} $F_\epsilon: Y \f X$. 

Moreover, if $g$ and $g_0$ are normalized so that
$
\Entv\left( Y,\, g\right) =\Entv\ \left( X,\, g_{0}\right),
$
then $F_\epsilon$ is almost-isometric. Furthermore if $\epsilon=0$, then $F_0$ is an isometric covering.  
\end{theorC}

The aim of the present paper is to analyze to what extent the analogous of Theorem C holds true in the \K\ setting by substituting  the volume entropy with the \emph{diastatic entropy}. 

In order to state Theorem \ref{thm dentropy2}  we need the following definitions analogous to those needed in the statement of Theorem C.
We say that  a \K\  manifold  $\left(Y,\, g\right)$ has {\em bounded diastatic  Hessian} if, for any point  $p$
of  its   universal \K\  covering $\left(\w Y,\, \w g\right)$ the following two conditions hold true:
\begin{equation}\label{ex(a)}
\int_{\w Y}  \w \rho_p \left( q \right) e^{-c\, \D_p\left( q\right) }\, \nu_{\w g}(q)< \infty, \qquad \forall \,c>\frac{\Ent\left( Y,\, g\right)}{\X \left( \w g \right)};
\end{equation}
there exists a positive constant $C$ such that   
\begin{equation}\label{ex(b)}
\left|\lmb^p_j\right| < C,\ \ j=1,\dots,m,
\end{equation} 
where $\lmb^p_j$ are the eingenvalues
the Hessian of the diastasis $\D_p$.
\begin{thm}\label{thm dentropy2} 
Let $\left(Y,\, g\right)$ be a compact \K\ manifold of  complex dimension $n\geq 2$
with bounded diastatic Hessian 
 and let  $\left(X,\, g_0\right)$  be a  compact complex hyperbolic manifold of the same dimension of $Y$.
 If $f:Y \f X$ is a non zero degree continuous map  and  there exists   a sufficiently small positive constant $\epsilon$ such that
\begin{equation}\label{EQ perto}
\Ent^{2n}\left( Y,\, g\right)\, \Vol\left( Y,\, g\right) - \left| \deg \left( f \right) \right| \, \Ent^{2n} \left( X,\, g_{0}\right)\, \Vol \left( X,\, g_{0}\right)  \leq\epsilon,
\end{equation}
then $f$ is homotopic to a $\mathcal {C}^1$-{covering} $F_\epsilon: Y \f X$. 
Moreover, if $g$ and $g_0$ are normalized so that
$
\Ent\left( Y,\, g\right) =\Ent \left( X,\, g_{0}\right),
$
then $F_\epsilon$ is almost-isometric. Furthermore if $\epsilon=0$, then $F_0$ is a holomorphic or anti-holomorphic isometric covering.   
\end{thm}

\begin{rmk}\rm
We believe that the map $F_\epsilon$ in Theorem \ref{thm dentropy2} is indeed a diffeomorphism and that condition \eqref{ex(a)} is  redundant.
\end{rmk}

Conditions \eqref{ex(a)} and \eqref{ex(b)} are somehow technical, so it is natural to seek for more topological and geometrical ones 
yielding {to} the same conclusions of Theorem \ref{thm dentropy2}.
This is achieved in Theorem \ref{thm proper} below which represents our second result.
One first topological condition is the following.
Let  $(M_1, \, g_1)$ and  $(M_2, \, g_2)$ be two  Riemannian manifolds.
We will say that   $M_1$ is a \emph {strongly proper} submanifold of $M_2$ if there exists  an isometric immersion
$\varphi: (M_1, \, g_1) \f (M_2, \, g_2)$, called a {\em strongly proper map}, such that one of its lift $\w \varphi: (\w M_1, \, \w g_1) \f (\w M_2, \, \w g_2)$  to the Riemannian universal covering manifolds 
satisfies the following condition: 
 {\em for any $\mu > 0 $ and $q\in \w M_1$, there exist two constants $L_1$ and $L_2$, such that
\begin{equation}\label{strongly proper}
\tilde\rho_1 \left(q,\, p\right) < L_1 \, e^{\,\mu\, \tilde\rho_2 \left(\w \varphi\left(q\right),\, \w \varphi\left(p\right)\right)} + L_2, \qquad \forall \, p \in \w M_{1}, 
\end{equation}
where $\tilde\rho_1$ and $\tilde\rho_2$ are the geodesics distances on $(\w M_1, \, \w g_1)$ and $(\w M_2, \, \w g_2)$ respectively.}
Notice that the previous  definition does not depend on the chosen lift and that an isometric immersion
$\varphi: (M_1, \, g_1) \f (M_2, \, g_2)$ is strongly proper  if there exists a polynomial $P$ such that  
$\rho_1 \left(q,\, p\right) <  P \left( \rho_2 \left(\w\varphi\left(q\right),\,\w \varphi\left(p\right)\right)\right).$ 
\begin{thm}\label{thm proper} 
Let $\left(Y,\, g\right)$ be a compact \K\ manifold of  complex dimension $m\geq 2$
which is a strongly proper  \K\  submanifold of a classical local hermitian symmetric space of non compact type and
let  $\left(X,\, g_0\right)$  be a  compact complex hyperbolic manifold of the same dimension of $Y$.
If $f:Y \f X$ is a non zero degree continuous map  satisfying \eqref{EQ perto} above, then the same conclusions of Theorem \ref{thm dentropy2} holds true.
\end{thm}

The paper is organized as follows. In section 2 after recalling the main properties of Calabi's diastasis function and diastatic hessian, we focus on the properties of hermitian symmetric spaces of noncompact type needed in the proof of the main results. Section 3 is dedicated to the definition and main properties of the barycentre map in the \K\ setting. Finally Section 4 contains the proof of Theorem 1 and 2.

\vskip 0.3cm

\noindent{\bf Acknowledgments.} The author would like to thank Professor Sylvestre Gallot and Professor Fabio Zuddas for various stimulating discussions and their valuable comments. The author gratefully thanks the referee for the constructive comments and recommendations which definitely help to improve the readability and quality of the paper. 

\section{Diastatic hessian and HSSNCT}
First briefly recall the definition of diastasis function.
Let $\left( \w Y,\, \w g \right)$ be a real analytic \K\ manifold, namely a complex manifold $\w Y$ endowed with a real analytic \K\ metric $\w g$. A real analytic \K\ metric $\w g$ is characterized for the local existence of a real analytic function $\Phi:V \f \R$, called   \K\ potential, such that $\w  \omega_{|V} =\frac{i}{2}\, \partial\ov\partial\, \Phi$, where $\w \omega$ is the \K\ form associated to $\w g$. Let $z=(z_1,\dots,z_n)$ be local coordinates  around a point $p_0\in V$, by duplicating the variables $z$ and $\ov z$
the real analytic \K\ potential $\Phi$ can be complex analytically
continued to a function $\hat \Phi: U\times U \f \C$ defined  in a neighborhood
$U\times U \subset V \times V$ of 
$\left(p_0, p_0\right)$ which is holomorphic in the first entry and antiholomorphic in the second entry. 
E. Calabi in its seminal paper \cite{calabi}, introduced the {diastasis function} $\D:U \times U \f \R$, the \K\ invariant defined by:
\begin{equation}\label{diast func}
\D\left( p,\, q \right):=\hat \Phi \left({z(p)},\, \ov {z(p)}\right)+\hat \Phi \left({z(q)},\,  \ov {z(q)}\right)-\hat \Phi \left({z(p)},\,  \ov {z(q)}\right)-\hat \Phi \left({z(q)},\,  \ov {z(p)}\right).
\end{equation}
One can see that it is uniquely determined by the \K\ metric $\w g$, i.e. does not depend on the choice of the \K\ potential $\Phi$ or on the local system of coordinates. Moreover, when we fix one of its entries, let's say $p$, then the \emph{diastasis centred in $p$}, $\D_p: U \f \R$ given by
$
\D_p\left( q \right):=\D\left(p,\, q \right)
$
is a \K\ potential. The reader is referred to \cite{LoiZedda-book} for further details and for  an updated account on projectively induced \K\ metrics.

In the proof of our results we need the following two lemmas about the diastasis function and Proposition \ref{prop hypo thm} that summarize the properties of classical Hermitian symmetric spaces of non compact type (from now on HSSNCT). {The interested reader can find in \cite{mossau} and \cite{mossaent} a computation of the diastatic entropy and the volume entropy of a HSSNCT}.
\begin{lem}[E. Calabi \cite{calabi}]\label{hereditary00000}
Let $\psi: (M_1,\, g_1) \f (M_2, \, g_2)$ be a holomorphic and isometric immersion between \K\ manifolds and suppose that $g_2$ is real analytic. Then $g_1$ is real analytic and for every couple of points $p,q \in M_1$ 
\[
\D^{M_1}\left( p,\, q\right)  =\D^{M_2}\left( \psi \left( p\right) ,\, \psi \left( q\right) \right) ,
\]
where $\D^{M_1}$ and $\D^{M_2}$ are respectively the diastasis of $(M_1,\, g_1)$ and $(M_2, \, g_2)$.
\end{lem}
\begin{lem}\label{cor hered}
Let $\psi: (M_1,\, g_1) \f (M_2, \, g_2)$ be a holomorphic and isometric immersion between \K\ manifolds and suppose that $M_2$ has globally defined diastasis $\D^{M_2}: M_2 \times M_2 \f \R$. Then $M_1$ has globally defined diastasis $\D^{M_1}:M_1 \times M_1 \f \R$ given by
\begin{equation}\label{hered glob diast}
\D^{M_1}\left( p,\, q\right)  =\D^{M_2}\left( \psi \left( p\right) ,\, \psi \left( q\right) \right).
\end{equation}
In particular the gradients and the hessians of $\D^{M_1}$ and $\D^{M_2}$ are {(locally)} related by the following identities:
\begin{equation}\label{GradDiast}
\psi_* \left( \grad_p \D^{M_1}_q \right) = 	\pi \left( \grad_{\psi (p)} \D^{M_2}_{\psi(q)}  \right)
\end{equation}
where $\pi: T _{\psi (p)} M_2 \f \psi _* \left(T_pM_1\right)$ is the orthogonal projection, and
\begin{equation}\label{HessDiast}
\Hess_{\psi(p)}\D^{M_2}_{\psi(q)}\left(\psi_*\xi,\, \psi_*\zeta\right) =  \Hess_p \D^{M_1}_q (\xi, \, \zeta) + \operatorname{II}_p\left( \xi,\, \zeta \right)\D^{M_2}_q,
\end{equation}
where $\operatorname{II}_{p}$ is the second fundamental form at $p \in M_{1}$.
\end{lem}
\proof
Equality  \eqref{hered glob diast} is an immediate consequence of Lemma \ref{hereditary00000}. 
Equality \eqref{GradDiast} is easily achieved: let  $E_1,\dots, E_{2n}$  be an orthonormal basis of $T_pM_1$, $p\in M_1$
then 
\[
\psi_* \left( \grad_p \D^{M_1}_q \right)= \sum_{k=1} ^{2n} \left( E_k\, \D^{M_1}_q \right)\psi_*E_k
= \sum_{k=1} ^{2n} \left( E_k\, \D^{M_2}_{\psi (q)}\circ \psi \right)\psi_*E_k
\]
\[
= \sum_{k=1} ^{2n} \left( \psi_* E_k \right) \D^{M_2}_{\psi (q)}\psi_*E_k = 	\pi \left( \grad_{\psi (p)} \D^{M_2}_{\psi(q)}  \right).
\]
It remains to prove \eqref{HessDiast}
For any $\xi, \, \zeta \in T_p M_1$ we have
$$
\Hess_p\D^{M_1}_q(\xi, \, \zeta) = \xi\left(  \zeta\,  \D^{M_1}_q\right) - \left(\nabla^{M_1}_\xi\zeta \right)\D^{M_1}_{q}
$$
and
$$
\Hess_{\psi(p)}\D^{M_2}_{\psi(q)}(\psi_*\xi,\, \psi_*\zeta) = \psi_* \xi\left(\psi_*  \zeta\,  \D^{M_2}_{\psi(q)}\right) - \left(\nabla^{M_2}_{\psi_* \xi} {\psi_* \zeta}\right) \D^{M_2}_{\psi(q)}
$$
hence
$$
\Hess_{\psi(p)}\D^{M_2}_{\psi(q)}(\psi_* \xi, \, \psi_*  \zeta) - \Hess_p\D^{M_1}_q(\xi, \, \zeta) = \left(\nabla^{M_2}_{\psi_* \xi} \psi_* \zeta - \psi_* \nabla^{M_1}_\xi \zeta  \right) \D^{M_2}_{\psi(q)}
$$
$$
=  \operatorname{II}_q\left( \xi,\, \zeta \right)\D^{M_2}_{\psi(q)}.
$$
\endproof
\begin{prop}\label{prop hypo thm}
Let $\left(\Omega, g^{\Omega}\right)$ be a HSSNCT, with $g^{\Omega}$ normalized in order to have holomorphic sectional curvature between $0$ e $-4$, then 
\begin{itemize}
\item the diastasis $\D^\Omega$ and the geodesic distance $\rho^\Omega$ are related by the following inequality
 \begin{equation}\label{diast poly e distantance}
\D^\Omega(w,z)
\geq 2 \log\cosh \left( \rho^\Omega(w,\, z)\right);
\end{equation}
\item if $(\Omega, g^{\Omega})$ is of classical type, then
\begin{equation}\label{grad class thm}
\mathcal X \left({ g^{\Omega}} \right)= \sup_{p,\, q\,\in\, \Omega}{\|\grad_p \D^{\Omega}_q\|}<\infty.
\end{equation}
Moreover the eigenvalues of the hessian of the diastasis are bounded, more  precisely for any  ${z, \, y \in \Omega}$ and any unitary  $\xi\in T_z\Omega$, we have 
\begin{equation}\label{hessian bounded}
0 < \Hess_z \D_y ^\Omega\left(\xi,\,\xi \right)<4.
\end{equation}
\end{itemize}
\end{prop}
\proof
We firstly consider the case of a HSSNCT of rank one, namely the complex hyperbolic space $(\C H^n,\,\w g_{\hyp} )$.
Let $\C H^n=\left\{ z \in \C^n : \|z\|^2=|z_1|^2+\dots+|z_n|^2<1 \right\}$ be the unitary disc endowed with the hyperbolic metric $\w g_{\hyp}$ of constant holomorphic sectional curvature $-4$. The associated \K\ form is $\w\omega_{\hyp}=-\frac{i}{2}\,\de\deb\,\log\left( 1-\|z\|^2\right)$ and the diastasis is given by
\begin{equation}\label{diast hyp}
\D^h(w,z)=-\log\left(\frac{\left(1-\|z\|^2\right)\left(1-\|w\|^2\right)}{\left|1-z w^*\right|^2}\right).
\end{equation}
Recalling the expression of geodesic distance,
\begin{equation}\label{geodesic distance hyperb}
\w \rho_{\hyp}\left(w ,\, z\right) = {\arctanh\left( \left|{\frac{w-z}{1-zw^*}}\right|\right) }
\end{equation}
we can conclude that the distance and the diastasis of the complex hyperbolic space are related by
\begin{equation}\label{diast and metric hyp}
\D^h(w,z)=2\log\cosh \left(\w \rho_{\hyp}(w,z)\right).
\end{equation}

By the polydisc theorem (see e.g. \cite{helgason}), for any couple of points $p, q \in \Omega$ there exists a totally geodesic polydisc $\left(P, g^P\right)$ of dimension $r = \rank \Omega$, holomorphically imbedded in $\Omega$ such that $p,q \in P$. By a  $r$-dimensional \emph{polydisc} $\left(P, g^P\right)$ we mean the following product of one dimensional  complex hyperbolic spaces with holomorphic sectional curvature $-4$,
\begin{equation}\label{poly defn} 
\left(P, g\right)=\left(\C H^1, \w g_h \right) \times \dots \times \left(\C H^1, \w g_h \right),
\end{equation}
where $P=\left\{(z_1, \dots, z_r) \in \C^r \, : \, \left|z_j\right| < 1, \ j=1,\dots,r\right\}$. The diastasis is the sum of the diastasis of each factor:
\begin{equation}\label{diast poly}
\D^P(w,z)=- \sum_{j=1}^r \log\left(\frac{\left(1-|z_j|^2\right)\left(1-|w_j |^2\right)}{\left|1-z_j  \ov w_j\right|^2}\right).
\end{equation}
By \eqref{geodesic distance hyperb} we see that the geodesic distance of $P$ is given by
\begin{equation}\label{geodesic distance poly}
\rho^P \left(w ,\, z\right) =\sqrt{ \sum_{j=1}^r \w \rho^{2}_{\hyp}\left(w_j ,\, z_j\right)} =\sqrt{ \sum_{j=1}^r  {\arctanh^2\left( \left|{\frac{w_j-z_j}{1-z_j\ov w_j}}\right|\right) }}.
\end{equation}
Using \eqref{diast and metric hyp} we obtain the following inequality
\begin{equation*}%\label{diast poly e distantance}
\begin{array}{C}
\D^P(w,z)=\sum_{j=1}^r \D^h(w_j,z_j)= 2 \sum_{j=1}^r \log\cosh \left(\w \rho_{\hyp}(w_j,\, z_j)\right) \\[1em]
\geq 2 \log\cosh \left( \rho^P(w,\, z)\right).
\end{array}
\end{equation*}
Inequality \eqref{diast poly e distantance} follows by combining  the previous inequality, the polydisc theorem, Lemma \ref{hereditary00000} and the fact that a HSSNCT has  globally defined diastasis (see for example \cite{loi hss}).

We first  prove \eqref{grad class thm} and \eqref{hessian bounded} for  the first classical domain  
$$\Omega_1=\Omega_1[m,m]=\left\{ Z \in M_{m,m} \, : \, \det \left(I - ZZ^*\right) >0  \right\}$$ 
endowed with its symmetric metric $g^{\Omega_1}$ of holomorphic sectional curvature between $0$ and $-4$. The \K\ form associated to $g^{\Omega_1}$ is
$\omega^{\Omega_1} = - \frac{i}{2} \de \deb \log \det \left( I - ZZ^* \right)$.
The diastasis centered in the origin is given by 
\begin{equation}\label{diast first dom}
\D_0^{\Omega_1}\left( Z \right) = - \log \det \left( I - ZZ^* \right).
\end{equation}
A straightforward computation show that
\begin{equation}\label{eq diff diag diast}
\begin{array}{C}
d_Z\D_0^{\Omega_1} =\sum_{h,k=1}^m\left(\Tr \left[\left(I-ZZ^*\right)^{-1}Z \, E_{kh}\right] d \ov z_{hk} \right.
\\[1em]
\left. + \Tr \left[\left(I-ZZ^*\right)^{-1}  E_{hk} \, Z \right] d z_{hk}\right)
\end{array}
\end{equation}
and
\begin{equation}\label{eq form first mm}
\begin{array}{C}
\omega^{\Omega_1} =\frac{i}{2}\sum_{i,j,h,k=1}^m\Tr \left[ \left(I-ZZ^*\right)^{-1}E_{ij}\, Z^* \left(I-ZZ^*\right)^{-1} Z\, E_{kh}  \right.
\\[1em]
\left. + \left(I-ZZ^*\right)^{-1} E_{ij}\, E_{kh} \right] dz_{ij} \wedge d \ov z_{hk},
\end{array}
\end{equation}
where $z_{11}, z_{12}, \dots, z_{mm}$ are the standard coordinates of $M_{m,m}$ denoting the entries of the matrix $Z$ and $E_{kh}$ is the matrix  with all the entries zero but the $kh$-th equal to one.

Since the group of holomorphic isometries $\mathcal G =\Isom\left({\Omega_1}, \, g^{\Omega_1}\right) \cap \Aut\left(\Omega_1\right)$ acts transitively on $\Omega_1$, by Lemma \ref{cor hered}, we can study $\grad \D^{\Omega_1}_W$ and $\Hess_Z \D^{\Omega_1}_W$, assuming $W=0$. Moreover, given unitary matrices $U_1,U_2 \in \U(m)$ the map $Z \mapsto U_1ZU_2$ is a holomorphic isometry of $\left({\Omega_1}, \, g^{\Omega_1}\right) $, that fixes the origin. Let $P'$ be the totally geodesic \K\ embedded $m$-dimensional polydisc of equation  $ P'= \left\{ Z \in \Omega_1 \, : \, z_{ij}=0 \text{ if } i\neq j\right\}$ (notice that $m$ is the rank of $\Omega_1$). Since $U_1,\,U_2$ can be choosed so that $V=U_1ZU_2$ is diagonal, by applying once again Lemma \ref{cor hered}, we can assume $Z \in P'$.

A straightforward computation shows that the gradient and the hessian of the diastasis restricted to $P'$ are given respectively by:
\begin{equation}\label{grad first diag}
\left.\grad \D^{\Omega_1}_0\right|_{P'}=2 \sum_{j=1}^m \left(1-\left|z_{jj}\right|^2  \right)\left( z_{jj} \frac{\de}{\de z_{jj}} + \ov z_{jj} \frac{\de}{\de \ov z_{jj}}\right)
\end{equation}
and
\begin{equation}\label{hessian polyy}
\begin{array}{C}
\left.{\Hess_Z \D^{\Omega_1}_0}\right|_{P'}= \sum_{j,k=1}^m \left(  \frac{dz_{jk}\otimes d\ov z_{jk}+d\ov z_{jk}  \otimes d z_{jk}}{\left(1-|z_{kk}|^2\right)\left(1-|z_{jj}|^2\right)}\right.  \\[1em]
\left.-\frac{\ov z_{jj} \, \ov z_{kk}\, dz_{jk}\otimes d  z_{jk} + z_{jj} \,  z_{kk} \,d\ov z_{jk}\otimes d \ov z_{jk}}{\left( 1-| z_{jj}|^2\right)\left( 1-| z_{kk}|^2\right)}\,\delta_{jk}\right).\\[1em]
\end{array}
\end{equation}
By the previous argument we can suppose $Z \in P'_+=\{Z \in P'  \, | \, z_{jj}\geq0, \, j=1,\dots,n\}$ and easily conclude that
\begin{equation}\label{grad class prof}
\mathcal X \left({\w g^{\Omega_1}} \right)= \sup_{p,\, q\,\in\, \w Y}{\|\grad_p \D^{\Omega_1}_q\|}<2 \sqrt n.
\end{equation}
Consider the orthonormal basis of $T_Z \Omega_1$, $$u_{jk}=\sqrt{\left( 1-| z_{jj}|^2\right)\left( 1-| z_{kk}|^2\right)}\left( \frac{\de}{\de z_{jk}} + \frac{\de}{\de \bar z_{jk}}\right)$$ and $u_{n+j n+k}=Ju_{jk},$ $j,k=1,\dots, n$ and notice that ${\nabla d_Z \D^{\Omega_1}_0}_{|_{P'_+}}\left( u_{jk}, \, u_{ls}\right)$ is a diagonal matrix with eigenvalues $0<\lmb_{jk}<4$.
Thus,  we conclude that for $Z,W \in {\Omega_1}$ and any unitary $\xi \in T_Z \Omega_1$
%\begin{equation}\label{norm grad first}
%\sup_{Z,W \in {\Omega_1} }\left \| \grad _Z\D_W \right\|^2 =  4 \, m 
%\end{equation}
%and
\begin{equation}\label{norm hess first}
0 < \nabla d_Z \D^{\Omega_1}_W \left(\xi, \xi \right)<4.
\end{equation}

We can address now the general case. Let $(\Omega, g^{\Omega})$ be any classical HSSNCT. It is known that $(\Omega, g^{\Omega})$ can be complex and totally geodesic embedded into $\Omega_1[m, m]$, for $m$ sufficiently large (this is obvious for the domains $\Omega_1$, $\Omega_2$ and $\Omega_3$, while for the domain $\Omega_4$, associated to the so called Spin-factor, the explicit embedding can be found at the bottom of p. 47 in \cite{crimmon}). Hence  by Lemma \ref{cor hered}, \eqref{grad class prof} and \eqref{norm hess first} we deduce the validity of \eqref{grad class thm} and \eqref{hessian bounded}. The proof of Proposition \ref{prop hypo thm} is complete. 
\endproof

\begin{cor}\label{lem hessian of the diast hyperb} 
Let $\left( \C H^n, \, \w g_{\hyp} \right)$ be the complex hyperbolic space  with associated diastasis $\D ^h$  (see \eqref{diast hyp}). Denoted by $J$ the complex structure, the Hessian $\nabla d \D^h $ of the diastasis can be written
\begin{equation}\label{hessian compl hyperb}
\begin{array}{C}
\nabla d_z\D^h_w =\\ [1ex]
= 2\,  {{\w g}_{\hyp}}\left( z\right) - \frac{1}{2} d_z \D^h_w
 \otimes d_z \D^h_w +\frac{1}{2}(d_z\D^h_w \circ J_z)\otimes(d_z\D^h_w \circ J_z)),\end{array}
\end{equation}
for all $z,w \in \C H^n$.
\end{cor}
\proof Consider $\left( \C H^n, \, \w g_{\hyp} \right)$ realized as the holomorphic and totally geodesic submanifold of $\Omega_1[n,n]$ of equation $z_{jk}=0$ if $j>1$. Observe that the diastasis centered in the origin of $\left( \C H^n, \, \w g_{\hyp} \right)$ is the restriction of $\eqref{diast first dom}$ to $\C H^n=\{Z \in \Omega_1[n,n]\ : \ z_{jk}=0, \ \forall \, j\neq 1 \}$, i.e.
$
\D^h_0=- \log \det \left( 1 - \sum_{j=1}^n \left| z_{1j} \right|^2 \right).
$

Notice that the group of holomorphic isometries of $\left( \C H^n, \, \w g_{\hyp} \right)$ acts transitively on $\C H^n$ and that it contains $\U(n)$. Therefore, in order to prove \eqref{hessian compl hyperb}, arguing as above we see that  it is enough to assume $w=0$ and $z$ with $z_{12}=\dots=z_{1n}=0$. By \eqref{eq diff diag diast}, \eqref{eq form first mm}  and \eqref{hessian polyy}, we see that

\begin{equation*}%\label{hessian compl hyp}
\begin{array}{C}
{\nabla d_z \D^h_0}=  \frac{\sum_{j,k=1}^n \left( dz_{jk}\otimes d\ov z_{jk}+d\ov z_{jk}  \otimes d z_{jk} \right) - \ov z_{11}^2 \, dz_{11}\otimes d  z_{11} - z_{11}^2  \,d\ov z_{11}\otimes d \ov z_{11}}{\left(1-|z_{11}|^2\right)^2}\\[1em]
=2\,  {{\w g}_{\hyp}}\left( z\right) - \frac{1}{2} d_z \D^h_0
 \otimes d_z \D^h_0 +\frac{1}{2}(d_z\D^h_0 \circ J_z)\otimes(d_z\D^h_0 \circ J_z)).
\end{array}
\end{equation*}
\endproof

\section{The barycentre map $\w F_c$}\label{sec bary}
Let  $\left( Y,\, g\right) $ be a  compact \K\ manifold with  universal \K\ covering $\left( \w Y, \, \w g\right)$ having globally defined diastasis. We define a positive finite measure $d\mu^c_y$ on $\w Y$ by
\begin{equation}\label{family of metrics}
d\mu^c_y\left( z\right) =e^{-c\, \D_y\left( z\right) }\, \nu_{\w g}, \quad c>\frac{\Ent\left( Y,\, g\right)}{\X \left( \w g \right)}.
\end{equation}
Let $(X, g_0)$ be a compact complex hyperbolic manifold of the same dimension of $Y$,  $f:Y \f X$ be a continuous map and let $\w f:\w Y \f \C H^n$ be its lift to the universal covers.
\begin{defn}\label{bary map} \rm For any  $c >\frac{\Ent\left( Y,\, g\right)}{\X \left( \w g \right)}$, we define the \emph{barycentre map} $\w F_c : \w Y \f \C H^n$, as the map that associates at $y \in \w Y$ the point where the function $\B_y:\C H^n\f\R^+$
\begin{equation}\label{bary}
x \mapsto \int_{\w Y} \D^h\left(\w f(z),\, x\right)\, d\mu^c_y(z)
\end{equation}
attains its unique point of minimum. 
\end{defn}
Here the notion of barycentre used by G. Besson, G. Courtois and S. Gallot in  \cite{bcg1} has been modified using in \eqref{bary} the Calabi's diastasis function $\D^h$ instead of the distance ${\w\di}_h$.
The following result  assures us that  the barycentre map $\w F_c$ is indeed  well defined.
\begin{lem}\label{prop by}
The function $\B_y:\C H^n\f\R^+$ admits a unique point of minimum.
\end{lem}
\proof
First we need to prove that $\B_y$ is well defined, namely that \eqref{bary} is convergent. Since $X$ and $Y$ are compact, by standard Riemannian geometry we can prove that, for given $x\in X$ and $y \in Y$, there exist constants $C_1$ and $C_2$ such that
$
\w \rho_h \left( x, \, \w f \left( z \right) \right)  \leq C_1\,\w \rho \left( y, \, z \right)  + C_2.
$
Therefore, for $\w \rho \left( y, \, z \right) >> 0$ there exists a positive constant $C_3$, such that:
\[
\D^h\left( x, \, \w f \left( z \right) \right) = 2\, \log \cosh \left( \w \rho_h \left( x, \, \w f \left( z \right) \right) \right) \leq 2\, \log \cosh \left( C_1\,\w \rho \left( y, \, z \right)  + C_2 \right) 
\]
\[
 \leq  C_3 \, \w \rho \left( y, \, z \right),
\]
where in the first equality we use \eqref{diast and metric hyp} and in the last inequality the fact that $\lim_{t \f +\infty}\frac{\log\cosh t} { t} = 1$. 
By \eqref{ex(a)}, we conclude that
\[
 \int_{\w Y} \D^h\left(\w f\left( z\right) ,\,x\right)\, d\mu^c_y\left( z\right) < \int_{\w Y}  C_3 \, \w \rho \left( y, \, z \right) \, d\mu^c_y\left( z\right)  < \infty,
\]
i.e. \eqref{bary} is well defined. 

We show now that the function $\B_y$ admits a point of minimum. Since $\| \grad_z \D^h \| = 2 \,\| z \| < 2$ for any $z \in \C H^n$, by the theorem of derivation under integral sign, we have
\begin{equation*}%\label{grad By} 
\grad_x \B_y = \int_{\w Y} \grad_x \D^h_{\w f\left( z\right)}\, d\mu^c_y\left( z\right),
\end{equation*}
in particular, we see that $\B_y$ and $\grad_x \B_y$ are continuous.
Let $T$ be a bounded non empty open set of $\w Y$, and define
\[
K\left( x\right) =\min_{z\in T}\D^h\left( \w f\left( z\right) ,\, x\right),
\]
so
\[
\B_y\left( x\right)= \int_{\w Y} \D^h_{\w f\left( z\right) }\left( x\right) \, d\mu^c_y\left( z\right) \geq K\left( x\right) \int_{T} d\mu^c_y\left( z\right).
\]
By \eqref{diast and metric hyp} we see that $K\left( x\right) \f + \infty$ as $x \f \de\, \C H^n$, that is 
$ \B_y\left( x\right) \f  + \infty$ {  as  }$x \f \de\, \C H^n.$ Therefore $\B_y$ attains its minimum in $\C H^n$.

It remains to prove that the point of minimum is unique. Since $\w Y$ is a complete Riemannian manifold, it is enough to prove that $\B_y$
is a strictly convex function, that is, we have to prove that the hessian of 
$\B_y$ is positive definite. By \eqref{hessian bounded} we know that $\|\nabla d_z\D^h_w\| < \infty$ for any $z,$ $w \in \C H^n$, so by the theorem of derivation under integral sign, the hessian of 
$\B_y$ is continuous and given by
\[
\nabla d_z\B_y= \int_{\w Y}  \nabla d_z \D^h_{\w f\left( z\right) }\, d\mu^c_y\left( z\right).
\]
By \eqref{hessian bounded}, we see that $\nabla d_z \D^h_{\w f\left( z\right) }$ and $\nabla d_z\B_y$ are positive definite. The proof is complete. \endproof

The main properties of the barycentre map $\w F_c : \w Y \f \C H^n$ are described by the following proposition.

\begin{prop}\label{the diastatic barycentre map}
The barycentre map $\w F_c:\w Y \f \C H^n$
satisfies the following properties:
\begin{enumerate}
\item
it is a $\mathcal C^1$ map, characterized by the equation
\begin{equation}\label{implicit F}
d_{\w F_c(y)} \B_{y} = \int_{\w Y} d_{\w F_c(y)} \D^h_{\w f(z)}\, d\mu^c_y\left( z\right)=0;
\end{equation}
\item it is equivariant with respect to deck transformations and it descend to a $\mathcal C^1$ map 
\begin{equation}\label{def Fc}
F_c : Y \f X
\end{equation}
homotopic to $f: Y \f X$.
\end{enumerate}
\end{prop}
\proof
%%%%%%%%%%%%%%%%%%%%%%%%%%%%%%%%%%%%%%%%%%
%%%%%%%%%%%%%%%%%%%%%%%%%%%%%%%%%%%%%%%%%%
By Proposition \ref{prop by} it  follows  that $\w F_c \left( y\right)$ is characterized by the  equation 
\begin{equation*}%\label{implicit F}
d_{\w F_c\left( y\right)} \B_{y} = \int_{\w Y} d_{\w F_c(y)} \D^h_{\w f(z)}\, d\mu^c_y\left( z\right)=0.
\end{equation*}
In other terms, given an orthonormal basis $e_j$, we define the function $\Phi: \C H^n \times  \w Y \f \R^{2n}$ by 
$
\Phi \left( x, \,y\right)^j =  d_x \B_{y}(e_j).
$
  Then we have $\Phi \left( \w F_c \left( y \right) , \,y\right) =0$. Since $\mathcal X \left(g_0\right) <\infty$ and $\mathcal X \left(g\right)<\infty $ then $\|d_x \D^h_{\w f(z)}\, d_y  \D_z\| < \infty $ and by the theorem of derivation under the integral sign, the differential of $\Phi$ with respect to $y$ is given by 
\[
d_y \Phi \left( x, \,y\right) = - c \int_{\w Y} d_x \D^h_{\w f(z)}\, d_y  \D_z\, d\mu_y^c \left( z \right) \leq - c  \,  \mathcal X \left(g\right) d_x \B_y.
\]
Arguing as in the proof of Lemma \ref{prop by}, we see that the Hessian of $\D^h_{\w f(z)}(x)$ is bounded and positive definite and therefore the Jacobian matrix of $\Phi$ with respect to $x$ is continuous and positive definite at  $ \left( \w F_c \left( y \right) , \,y\right)$. Thus, we can apply the implicit function theorem and obtain the $\mathcal C^1$-regularity of the maps $F_c$.  This concludes the proof of \emph{(1)}.

Consider now $\Gamma = \pi_1 \left( Y,\, y_0\right)$ the group  of deck transformations of the universal covering of $Y$. The morphism
$ f_{\ast} :  \pi_1 \left( Y,\, y_0\right)  \to  \pi_1 \left( X,\, f\left(  y_0\right) \right) $ induces a representation 
$ r : \Gamma \to {\rm Isom} \left( \C H^n,\, \w g_0 \right) \cap \Aut\left( \C H^n\right) $ which satisfies 
$\w f \circ \gamma = r\left( \gamma\right)  \circ \w f$ for every $\gamma \in \Gamma$.
As $\gamma_{\ast} v_{\w g} = v_{\w g}$, and as $r\left( \gamma\right) $ is a holomorphic isometry of
$\C H^n$, we have, for every $y \in \w Y$ and every $x \in \w X$:
\[
\B_{\gamma\, y} \left(r\left( \gamma\right) x \right) =
\int_{\w Y} \D^h\left({\w f\left( z\right) },\, r\left( \gamma\right)  x\right) e^{-c\, \D\left( \gamma \, y,\, z\right) }\,\nu_{\w g}
\]
\[
=\int_{\w Y} \D^h\left({\w f\left( \gamma\, z\right) },\,  r\left( \gamma\right)  x\right) 
e^{-c \, \D\left( \gamma \, y,\,\gamma\, z\right) }\,\nu_{\w g}
\]
\[
= \int_{\w Y} \D^h\left( r\left( \gamma\right) \w f\left(z\right)  , \, r\left( \gamma\right)  x\right) e^{-c\, \D\left( y,\, z\right) }\,\nu_{\w g}
\]
\[
=\int_{\w Y} \D^h\left({\w f\left( z\right) },\, x\right) d\mu^c_y\left( z\right) =  \B_{y} \left( x \right)
\]
As $\B_{y}$ attains its minimum at the unique point $\w F\left( y\right) $,
this equality implies that $\B_{\gamma\, y}$ attains its minimum at the unique point 
$r\left( \gamma\right)  \w F\left( y\right) $. That is $\w F\left( \gamma \, y\right) = r\left( \gamma\right) \, \w F\left( y\right) $. Therefore $\w F_{c}$ is invariant with respect to deck transformations and it descends to a map
\begin{equation*}%\label{def Fc}
F_c : Y \f X.
\end{equation*}
In order to prove that the maps $F_c$ and $f$ are homotopic,
consider the Dirac measure $\delta_y\left( z\right) $ on $\w Y$.
Let us define the positive finite measure $d\mu_y^{c,\,t}$ as follows 
\[
d\mu_y^{c,\, t}\left( z\right) =t\, d\mu^c_y\left( z\right) +\left( 1-t\right) \delta_y\left( z\right) 
\]
and let  $\w F_{c,\, t} : \w Y \f \C H^n$ be the map given by 
\[
\w F_{c,\,t}\left( y\right) =\Bary\left(\w f_* d\mu^{c,\,t}_y\left( z\right) \right),
\]
i.e. $\w F_{c,\,t}\left( y\right) $ is the unique point where the function $ \B_{y,\,t}: \w X\f\R^+$ defined by
\begin{equation}\label{homobary}
\begin{split}
\B_{y,\,t} \left( x\right)  &=\int_{\w Y} \D^h\left( \w f\left( z\right) ,\, x\right) d\mu^{c,\,t}_y\left( z\right) =\\
&=t\int_{\w Y} \D^h\left( \w f\left( z\right) , \,x\right) d\mu_y\left( z\right) +\left( 1-t\right) \, \D^h\left(  \w f\left( y\right) ,\, x\right)
\end{split}
\end{equation}
attains its minimum. Clearly $\w F_{c,1}=\w F_{c}$. Let $\phi \in \Isom\left( \C H^n, \, \w g_0\right)  \cap \Aut\left( \C H^n\right)  $ such that $\phi \left( x\right) =0$, then
\begin{equation*}%\label{implicit F}
\D^h\left( x,\, z\right) =\D^h\left( 0,\,\phi\left( z\right) \right) =-\log\left( 1-\left|\phi\left( z\right) \right|^2\right) 
\end{equation*}
therefore $\D^h\left( x,\, z\right) \geq 0$ and $\D^h\left( x,\, z\right) = 0$ if and only if $x=z$, so the function $\B_{y,\, 0}$ attains its unique minimum for $x=\w f\left( y\right) $, i.e. $\w F_{c,\, 0}\left(  y\right) =\w f\left( y\right) $.

Arguing as before, we conclude that $\w F_{c,\, t}\left( y\right) $ is a well defined $\mathcal C^1$ map, equivariant with respect to deck transformations. So $\w F_{c,\,t}\left( y\right) $ descends to a homotopy $F_{c,\, t}\left( y\right) $ between $F_c$ and $f$.
\endproof

%\medskip

\section{The proof of Theorem \ref{thm dentropy2} and Theorem \ref{thm proper}.} \label{proof1}

Let $f: Y \f X$ be the continuous  function given in the hypothesis of Theorem \ref{thm dentropy2} and let $\w F_c: \w Y \f X$ be the associated barycentre map, given by Definition \ref{bary map}.

In order to differentiate \eqref{implicit F} under the integral sign, note that by \eqref{hessian compl hyperb} and $\X(\w g_h)=2$, we get
\[
\max_{\left\| u \right\|= \left\| v \right\|=1}\left|\nabla\left(d_{\w F_c(y)} \D^h_{\w f(z)}\, e^{-c\, \D(y,z)}\right)\left( u, \, v \right) \right|
\]
\[
\leq \max_{\left\| u \right\|= \left\| v \right\|=1} \left( \left|\nabla   d \D^h_{\w f(z)} \left( d \w F_c \left(  u \right) , \, v \right)     \right| + \left| d \D^h\otimes d \D \left( u, \, v \right) \right|\right) e^{-c\, \D(y,z)} 
\]
\[
\leq \left(  6 \left\| d \w F_c \right\|+ 2\,\X \left( g\right) \right)  e^{-c\, \D(y,z)},
\]
by Proposition \ref{the diastatic barycentre map} the map $\w F_c$ descend to a map $F_c: Y \f X$, so, as $Y$ is compact, $\left\| d \w F_c \right\|$ is bounded. Hence the norm of the derivative of the integrand in \eqref{implicit F} is bounded by a constant function, which (by the hypothesis $c > \frac{\Ent\left( Y,\, g\right)}{\X \left( \w g \right)}$) is integrable. Thus, bsy standard measure theory, we can derive \eqref{implicit F} under the integral sign. For every $v \in T_{\w F_c(y)} \C H^n$ and $u \in T_y(\w Y)$, we get
\begin{equation}\label{eq main1}
\begin{array} {C}
\displaystyle \int_{\w Y} \nabla d_{\w F_c(y)} \D^h_{\w f(z)}(d_y \w F_c(u),v)\, d\mu^c_y(z)\\[1em]
\displaystyle \quad =c \int_{\w Y} d_{\w F_c(y)} \D^h_{\w f(z)}(v)\, d_y \D_z(u)\,d\mu^c_y(z).
\end{array}
\end{equation}

Let us denote by
$K$, $H$ and $H'$ the symmetric endomorphisms of $T_{\w F_c(y)} \C H^n$ and $T_y \w Y$ defined by
\[
\w g _ {\hyp} \left(K(v),w\right)=\frac{1}{\int_{\w Y}\,d\mu^c_y(z)}\int_{\w Y} \nabla d_{\w F_c(y)} \D^h_{\w f(z)}(v,w) \,d\mu^c_y(z),
\]
\[
\w g _ {\hyp} \left(H(v),w\right) = \frac{1}{\int_{\w Y}\,d\mu^c_y(z)}\int_{\w Y} d_{\w F_c(y)} \D^h_{\w f(z)}(v)\, d_{\w F_c(y)} \D^h_{\w f(z)}(w) \,d\mu^c_y(z),
\]
\[
\w g \left(H'(u),t\right) = \frac{1}{\int_{\w Y} \,d\mu^c_y(z)}\int_{\w Y} d_y \D_z(u)\,d_y \D_z(t) \,d\mu^c_y(z),
\]
where $v,w\in T_{F_c(y)} \C H^n$ and $u,t\in T_y \w Y$.\\
By the Cauchy-Schwarz inequality and \eqref{eq main1}, we deduce
\begin{equation}\label{ineqdet}
\left|\w g _ {\hyp} \left( K\circ d_y\w F_c\left( u\right) ,\, v\right) \right| \leq c\, \w g _{\hyp} \left( H\left( v\right) ,v\right) ^\frac{1}{2}\, \w g \left( H'\left(u\right) ,u\right) ^\frac{1}{2}
\end{equation}
\begin{lem}\label{lemdet}
With the previous notations we have
\begin{equation}\label{diseq1}
\left| \det K  \right| \left| \det (d_y \w F_c) \right|  \leq \left(\frac{\mathcal X^2 (g)\, c^2}{2\,n}\right)^{n}(\det H)^\frac{1}{2}
\end{equation}
and
\begin{equation}\label{eq HsuK}\begin{split}
\frac{\left(\det H\right)^\frac{1}{2}}{\det K}=\frac{\left(\det H\right)^\frac{1}{2}}{\det \left(2I-\frac{1}{2}H-\frac{1}{2}JHJ\right)}< \left(\frac{1}{2\,n}\right)^n.
\end{split}\end{equation}
\end{lem}
\proof
Let $\left\{v_j\right\}$ be an orthonormal basis
of $T_y \C H^n$ which diagonalizes the symmetric endomorphism $H$. Now, if $d_y\w F_c$ is not invertible, the inequality is trivial. Suppose that $d_y\w F_c$ has maximal rank. Let $u'_j=\left( K\circ d_y\w F_c\right) ^{-1}\left( v_j\right) $. By the Gram-Schmidt orthonormalization applied to $\left\{u'_j\right\}$, with respect the positive bilinear form $\w g \left( H'\left( \cdot \right) , \cdot \right)$, we get an orthogonal basis $\left\{u_j\right\}$ such that $\w g \left( u_j ,\, u_j \right)^{-\frac{1}{2}}=\lmb_j,$ $j=1, \dots, 2\,n$ are the eigenvalues of $H'$. Then
\begin{equation*}%\label{diseq1}
\left|\det K \right| \left|\det \left( d_y \w F_c\right) \right| = \prod_{j=1}^{2n} 
\left|\w g _ {\hyp} \left( K\circ d_y\w F_c\left( u_j\right) ,\,v_j\right) \right| \left( \det H' \right) ^{\frac{1}{2}}, 
\end{equation*}
hence, by \eqref{ineqdet}
\begin{equation}\label{eq lem 45}
\left|\det K \right| \left|\det \left( d_y \w F_c\right) \right|\leq c^{2n}  \left( \det H\right) ^\frac{1}{2} \left( \det H' \right) ^{\frac{1}{2}}
\end{equation}
\begin{equation*}
\leq c^{2n} \left( \det H\right) ^\frac{1}{2} \left(\frac{1}{2\,n} \, \tr H' \right)^{n}
\end{equation*}
\begin{equation*}
=\left(\frac{\mathcal X^2 (g)\, c^2}{2\,n}\right)^{n} \left( \det H\right) ^\frac{1}{2},
\end{equation*}
where we use that the eigenvalues of $H'$ are positive  and that for any orthonormal basis $\{e_1 , \ldots , e_{2n}\}$ of $T_y\w Y$
\[\sum_{i =1}^{2n}\ \w g \left( H'\left( e_i\right) ,\, e_i\right)  = \frac{1}{\int_{\w Y}\, d\mu^c_y\left( z\right) }\int_{\w Y} 
\left(\sum_{i =1}^{2n}\left(d_y \D_z\left( e_i\right) \right)^2\right)\, d\mu^c_y\left( z\right)  \leq \mathcal X^2 \left( g \right)  .\]
So \eqref{diseq1} is proved. 
By \eqref{hessian compl hyperb} we see that
%\begin{equation}\label{eq HsuK}\begin{split}
$
\frac{\left(\det H\right)^\frac{1}{2}}{\det K}=\frac{\left(\det H\right)^\frac{1}{2}}{\det \left(2I-\frac{1}{2}H-\frac{1}{2}JHJ\right)}.
$
%\end{split}\end{equation}
Consider the function $H \mapsto \frac{\left(\det H\right)^\frac{1}{2}}{\det \left(2I-\frac{1}{2}H-\frac{1}{2}JHJ\right)}$ defined over the group of symmetric matrices non negatively defined and with trace  $\leq 4$ and dimension $2n\times 2n$ with $n\geq 2$. By \cite{bcg1} Appendix B, attains its maximum at $H=\frac{2}{n}\,I$. Hence
$
%\begin{equation}\begin{split}\label{ineqHK}
\frac{(\det H)^\frac{1}{2}}{\det K}< \left(\frac{1}{2\,n}\right)^n.
%\end{split}\end{equation}
$
\endproof

 In order to prove Theorem \ref{thm dentropy2} notice that  the quantity $\Ent^{2n}\left(Y,\, g \right)\, {\Vol}\left(Y,\, g \right)$ is invariant by homotheties, hence it is not restrictive assume from the very beginning that $\Ent\left( Y,\, g\right) =\Ent \left( X,\, g_0\right)=\Ent \left( X,\, g_{\hyp}\right)=2\,n$. 
 The first part of Theorem \ref{thm dentropy2} will immediately follow by Theorem \ref{thm almost} below. 
{The second part of Theorem \ref{thm dentropy2} (the $\ep = 0$ case), is proved 
 in the last part of this section.}
  
\begin{thm}\label{thm almost}
Let $\left( Y,\, g\right)$ and $\left( X,\, g_{\hyp}\right)$ be as in Theorem \ref{thm dentropy2}.
Assume that $\Ent\left( Y,\, g\right) =\Ent \left( X,\, g_{\hyp}\right)$ and that
\begin{equation}\label{ipvol}
\Vol\left( Y,\, g\right)  < \left( 1+\ep\right) \left| \deg \left( f\right)  \right|\Vol\left( X,\, g_{\hyp}\right).
\end{equation}
If $\ep>0$ is small enough and $c$ is such that $\left(\left(\frac{\X(g)\,c}{2\,n}\right)^{2n}-1\right)< \frac{\ep}{(1+\ep)}$, then the map $F_c$ is a $\mathcal C^1$ covering map such that
\begin{equation}\label{eq quasi}
\begin{split}
A'\left( \ep\right)  \leq \frac{\left\|d_y F_c\left( u\right) \right\|}{\left\|u\right\|}\leq A''\left( \ep\right)  \qquad \forall\, y \in Y,\ \forall u \in T_y Y
\end{split}\end{equation}
where $A'\left( \ep\right) ,\, A''\left( \ep\right)  \f 1$ as $\ep \f 0$.
\end{thm}
In order to prove the theorem, we  need of the following five lemmata (Lemma \ref{lemma ye}-\ref{epineq02}).
\begin{lem}\label{lemma ye}
Let $Y_\ep=\left\{y\in Y:\left| \Jac F_c\left( y\right) \right| < \left( 1-\sqrt \ep\right)  \left( 1+\delta\right) \right\}$ where $\delta>0$ is defined by
\begin{equation}\label{defd}
\delta=\left(\frac{\X\left( g\right) \,c}{2\, n}\right)^{2n}-1.
\end{equation}
Then, for $\delta< \frac{\ep}{\left( 1+\ep\right) }$, we have
\begin{equation*}
\Vol \left( Y_\ep\right) <2\, \sqrt \ep \Vol\left( Y\right) .
\end{equation*}
\end{lem}
\proof
%By \eqref{ineqentr} we know that $\left| \Jac F(y)\right|\leq 1+\frac{\delta}{c_0}$,  
By \eqref{diseq1} and \eqref{eq HsuK} we know that $\left| \Jac F_c\right| < 1+ \delta$, by the definition of $Y_\ep$ we get
\begin{equation*}
\begin{split}
&\left( 1+\delta \right)\Vol \left(Y \setminus Y_\ep\right) + \left(1- \sqrt \ep \right) \left( 1+\delta \right)\Vol \left( Y_\ep\right)\geq \int_Y\left| \Jac F_c\right| \, \nu_{g} .
\end{split}
\end{equation*}
Using the hypothesi \eqref{ipvol} we obtain
\begin{equation*}
\begin{split}
\int_Y\left| \Jac F_c\right| \, \nu_{g}  \geq |\deg\left(  f\right) | \Vol \left( X\right)  > \frac{1}{1+\ep} \Vol\left( Y\right)  > \frac{ 1+\delta}{1+2\ep}\Vol\left( Y\right),
\end{split}
\end{equation*}
Where the last inequality follows by the assumption $\delta< \frac{\ep}{\left( 1+\ep\right) }$. 
Thus
\begin{equation*}
\begin{split}
&\Vol \left(Y \setminus Y_\ep\right) + \left(1- \sqrt \ep \right) \Vol \left( Y_\ep\right)> \frac{ 1}{1+2\ep}\Vol\left( Y\right) ,
\end{split}
\end{equation*}
and so
\[
\Vol \left( Y_\ep\right) < \frac{2\,\sqrt\ep}{1+2\,\ep} \Vol\left( Y\right)  <2\, \sqrt \ep \, \Vol\left( Y\right).
\]
As wished.
\endproof

Let us denote $\w Y_{\ep}=\pi^{-1}\left( Y_{\ep}\right) $. By the definition of $Y_{\ep}$ and  \eqref{eq lem 45}, we get
\begin{equation}\label{ineqdet3}
\begin{split}
&\left(1-\sqrt {\ep} \right) \left(1+\delta\right) \leq \left| \Jac \w F_c(y)\right|\leq \frac{c^{2n}\left(\det H\right)^\frac{1}{2}\left(\det H'\right)^\frac{1}{2}}{\left | \det K \right |}, \quad \forall\,y \in \w Y \setminus \w Y_{\ep}
\end{split}
\end{equation}
hence, by \eqref{eq HsuK} and \eqref{defd}, we deduce
\begin{equation*}%\label{ineqdet1}
\begin{split}
\left(\det H'\right)
\geq\left(\frac{\left(1-\sqrt {\ep}  \right) \left(1+\delta\right)}{c^{2n}}\, \frac{\left | \det K \right |}{\left(\det H\right)^\frac{1}{2}} \right)^2
\geq  \left(1-\sqrt {\ep} \right)^2\left(\frac{\X^2(g)}{2\,n}\right)^{2n}.
\end{split}
\end{equation*}
Since $\det H' \leq  \left( \frac{\tr H'}{2n} \right)^{2n} \leq \left( \frac{\X^2 (g)}{2\,n}\right)^{2n}$, we get
\begin{equation}\label{ineqdet2}
\begin{split}
\left(1-\sqrt {\ep} \right)^2\left(\frac{\X^2(g)}{2\,n}\right)^{2n}\leq\det H' \leq \left( \frac{\X^2 (g)}{2\,n}\right)^{2n}
\end{split}
\end{equation}
As the maximum of $H' \mapsto \det H'$ is obtained at $H'=\frac{\X^2 (g)}{2\,n}I$ by a principle of stability of the maximum (see \cite{bcg3} pag. 157), there exist a positive constant $B'(n)$ such that, for $\ep < \frac{1}{\left(2 B' (n)\right)^4}$
\begin{equation}\label{max1}
\left\| H' - \frac{\X^2 (g)}{2\,n}I \right\| \leq B'(n)\,\ep^{\frac{1}{4}}, \qquad \forall \, y \in \w Y \setminus \w Y_{\ep}.
\end{equation}
On the other hand by \eqref{ineqdet3} we obtain
\begin{equation*}%\label{ineqdet1}
\begin{split}
\frac{\det H^\frac{1}{2}}{\det K} \geq \frac{\left(1-\sqrt {\ep}\right)\left(1+\delta\right)}{c^{2n} \left(\det H'\right)^{\frac{1}{2}}}\geq \left(1-\sqrt {\ep}\right)\left( \frac{1}{2n}\right)^n,
\end{split}
\end{equation*}
Where the second inequality follows by \eqref{ineqdet2}. By \eqref{eq HsuK}  we get
\begin{equation*}%\label{ineqdet1}
\begin{split}
\left(1-\sqrt {\ep}\right)\left( \frac{1}{2n}\right)^n\leq\frac{\left(\det H\right)^\frac{1}{2}}{\det \left(2I-\frac{1}{2}H-\frac{1}{2}JHJ\right)}\leq\left(\frac{1}{2n}\right)^{n}.\end{split}
\end{equation*}
As we see before the maximum of $H\mapsto \frac{\left(\det H\right)^\frac{1}{2}}{\det \left(2I-\frac{1}{2}H-\frac{1}{2}JHJ\right)}$ is obtained for $H=\frac{2}{n}I$, so by a principle of stability of the maximum  (see \cite{bcg1}), there exist a positive constant $B''(n)$ such that, for $\ep < \frac{1}{\left(2 B'' (n)\right)^4}$, we have
\begin{equation}\label{max2}
\left\| H - \frac{2}{n}I \right\| \leq B''(n)\,\ep^{\frac{1}{4}}, \qquad \forall \, y \in \w Y \setminus \w Y_{\ep}.
\end{equation}
From now on, we benote $B(n)$ the maximum between $B''(n)$ and $B'(n)$.
\begin{lem}
If $\ep<\frac{1}{\left( {4\,B\left( n \right) }\right)^4 }$ and $c$ is such that $\delta=\left(\frac{\X(g)\,c}{2\,n}\right)^{2n}-1< \frac{\ep}{(1+\ep)}$ then,  $\forall \, y \in \w Y \setminus \w Y_{\ep}$, we have
\begin{equation}\label{epineq00}
\begin{split}
\left\| d_y \w F_c(u) \right\| \leq \frac{c\left(B\left( n\right) \ep^{\frac{1}{4}} + \frac{2}{n}\right)^{\frac{1}{2}}\left(B\left( n\right) \ep^{\frac{1}{4}} + \frac{\X^2 (g)}{2\,n}\right)^{\frac{1}{2}}\left\| u \right\|}{ \left(2 - B\left( n\right) \ep^{\frac{1}{4}}\right)}
\end{split}
\end{equation}
and
\begin{equation}\label{epineq01}
\begin{split}
\left\| d_y \w F_c\left( u\right)  \right\|\geq\left(
\frac{\left( \left( 1+\delta\right)\left( 1- \sqrt \ep \right)   \right)^{\frac{2\,n}{2\,n-1}}  \left(2 - B\left( n\right) \ep^{\frac{1}{4}}\right)}{c\left(B\left( n\right)  \ep^{\frac{1}{4}} + \frac{2}{n}\right)^{\frac{1}{2}}\left(B\left( n\right) \ep^{\frac{1}{4}} + \frac{\X^2 \left( g\right)}{2\,n}\right)^{\frac{1}{2}}}\right)^{2n-1}\left\|u\right\|
\end{split}
\end{equation}
\end{lem}
\proof
By \eqref{max2} we have
\begin{equation}\label{ineq k-2i}
\begin{array}{C}
\left\| K - 2I \right\| = \left\| \frac{1}{2} H - \frac{1}{2} J^{-1} H J \right\| \leq \left\| \frac{H}{2}  - \frac{I}{n} \right\| + \left\|   J^{-1}\left( \frac{I}{n} - \frac{H}{2} \right) J \right\| 
\leq B(n)\,\ep^{\frac{1}{4}}
\end{array}
\end{equation}
Note that
\begin{equation*}
\begin{array}{C}
\w g_{\hyp} (2\,v,\, w) - \w g_{\hyp}\left( K\, v,\, w\right)  \leq \left| \w g_{\hyp} \left( K\, v ,\, w\right)  - \w g_{\hyp}\left( 2\, v,\, w\right)  \right|\leq\\[1ex]
\quad\leq \left\| K - 2\,I \right\| \left\| v \right\| \left\| w \right\|
\leq B\left( n\right) \ep^{\frac{1}{4}}  \left\| v \right\| \left\| w \right\|
\end{array}
\end{equation*}
and so
\begin{equation*}
\begin{split}
\w g_{\hyp}\left( K\, v,\,w\right)  \geq \w g_{\hyp} \left( 2\,v,\,w\right)  - B\left( n\right) \ep^{\frac{1}{4}}  \left\| v \right\| \left\| w \right\|.
\end{split}
\end{equation*}
Setting $v=d_y \w F_c(u)$ and $w=\frac{d_y \w F_c(u)}{ \left\| d_y \w F_c(u) \right\| }$  we obtain
\begin{equation}\label{ineq dF01}
\begin{split}
\w g_{\hyp}\left(K \circ d_y \w F_c (u), \frac{d_y \w F_c(u)}{ \left\| d_y \w F_c(u) \right\| }   \right) \geq  \left\| d_y \w F_c(u) \right\|  \left(2 - B(n) \ep^{\frac{1}{4}}\right).
\end{split}
\end{equation}
By \eqref{max2}, we see that
\begin{equation*}
\begin{split}
\w g_{\hyp} \left( H (u), u \right) - \w g_{\hyp} \left( \frac{2}{n}\, u,u \right)\leq \left\| H - \frac{2}{n}I \right\| \left\|  u \right\|^2\leq B(n)\ep^{\frac{1}{4}}\left\|  u \right\|^2.
\end{split}
\end{equation*}
therefore
\begin{equation}\label{inhuu}
\begin{split}
\w g_{\hyp}\left(H(u),u\right) \leq  \left(B(n) \ep^{\frac{1}{4}} + \frac{2}{n}\right)\left\| u \right\|^2.
\end{split}
\end{equation}
On the other hand, by \eqref{max1}, we get
\begin{equation*}
\begin{split}
&\w g_{\hyp} \left( H' (w), w \right) - \w g_{\hyp} \left( \frac{\X^2(g) }{2\,n}\, w,w \right)\leq \left\| H - \frac{\X^2(g)}{2\,n}I \right\| \left\|  w \right\|^2\\
&\quad\leq B(n)\ep^{\frac{1}{4}}\left\|  w \right\|^2.
\end{split}
\end{equation*}
and so
\begin{equation}\label{inh'ww}
\begin{split}
\w g_{\hyp}\left(H'(w),w\right) \leq  \left(B(n) \ep^{\frac{1}{4}} + \frac{\X^2 (g)}{2\,n}\right)\|w\|^2.
\end{split}
\end{equation}
Substituting \eqref{ineq dF01}, \eqref{inhuu} and \eqref{inh'ww} in \eqref{ineqdet} we obtain

\begin{equation*}%\label{epineq}
\begin{split}
\left\| d_y \w F_c(u) \right\| \leq\frac{c\left(B(n) \ep^{\frac{1}{4}} + \frac{2}{n}\right)^{\frac{1}{2}}\left(B(n) \ep^{\frac{1}{4}} + \frac{\X^2 (g)}{2\,n}\right)^{\frac{1}{2}}\left\| u \right\|}{ \left(2 - B(n) \ep^{\frac{1}{4}}\right)} 
\end{split}
\end{equation*}
We proved equation \eqref{epineq00}.
Let $0<\left| \lmb_1\right|^2\leq \dots \leq \left| \lmb_{2n}\right|^2$ the eigenvalues of the symmetric endomorphism defined by $\left( d_y \w F_c \right)^t d_y \w F_c$. So
\[
0<|\lmb_1|\leq \dots \leq |\lmb_{2n}|\leq \frac{c\left(B\left( n\right)  \ep^{\frac{1}{4}} + \frac{2}{n}\right)^{\frac{1}{2}}\left(B\left( n\right)  \ep^{\frac{1}{4}} + \frac{\X^2 \left( g\right) }{2\,n}\right)^{\frac{1}{2}}}{ \left(2 - B\left( n\right)  \ep^{\frac{1}{4}}\right)} 
\]
moreover, by the definition of $Y_\ep$ follow that for every $y \in \w Y \setminus \w Y_\ep$ we have, $\prod_{j=1}^{2n} \left| \lmb_j\right| \geq \left( \left( 1+\delta\right)\left( 1- \sqrt \ep \right)   \right)^{{2\,n}}$, therefore
\[
|\lmb_1|\geq \frac{\prod_{j=1}^{2n} \left| \lmb_j\right|}{|\lmb_{2n}|^{2n-1}}\geq \left(
\frac{ \left( \left( 1+\delta\right)\left( 1- \sqrt \ep \right)   \right)^{\frac{2\,n}{2\,n-1}} \left(2 - B(n) \ep^{\frac{1}{4}}\right)}{c\left(B(n) \ep^{\frac{1}{4}} + \frac{2}{n}\right)^{\frac{1}{2}}\left(B(n) \ep^{\frac{1}{4}} + \frac{\X^2 (g)}{2\,n}\right)^{\frac{1}{2}}}\right)^{2n-1}
\]
we conclude that 
\begin{equation*}%\label{epineq}
\begin{split}
\left\| d_y \w F_c(u) \right\| \geq \left(
\frac{\left( \left( 1+\delta\right)\left( 1- \sqrt \ep \right)   \right)^{\frac{2\,n}{2\,n-1}} \left(2 - B(n) \ep^{\frac{1}{4}}\right)}{c\left(B(n) \ep^{\frac{1}{4}} + \frac{2}{n}\right)^{\frac{1}{2}}\left(B(n) \ep^{\frac{1}{4}} + \frac{\X^2 (g)}{2\,n}\right)^{\frac{1}{2}}}\right)^{2n-1}\|u\|,
\end{split}
\end{equation*}
we just proved  \eqref{epineq01}. The proof is complete.
\endproof

For every $y\in \w Y$, $u\in T_y \w Y$ and $v \in T_{\w F_c (y)} \w X$ we define
\begin{equation}\label{defnk'}
%\begin{split}
k'_y\left( u,\,v\right) =\frac{1}{\int_{\w Y}d\mu_y\left( z\right) }\int_{\w Y} d_{\w F_c\left( y\right) } \D^h_{\w f\left( z\right) }\left( v\right)  d_y \D_z\left( u\right)\,d\mu_y\left( z\right)  %\,e^{-c\, \D_z\left( y\right) }\frac{\w \omega^n\left( z\right) }{n!}
%\end{split}
\end{equation}
\begin{lem}\label{lemHessK}
There exist a universal constant $C$ such that
\begin{equation}\label{epineq0000}
\begin{split}
\left\| \nabla_w k'(u,v)\right\|\leq C \left\| u \right\|\left\| v \right\| \left( \left\| w \right\| + \left\| d \w F_c (w) \right\| \right). 
\end{split}
\end{equation}
\end{lem}
\proof
Assume for the moment that the following derivations under the integral sign are allowed, for every $w \in T_yY$ we have
\begin{equation}\label{epineq}
\begin{array}{C}
\nabla_w k'\left( u,\,v\right) \int_{\w Y}\ d\mu_y^c\left( z \right) \\[1em]
%\left( \int_{\w Y}\ d\mu_y^c(z) \right)\nabla_w k'(u,v)\\
\quad=\int_{\w Y} \nabla d_{\w F_c(y)} \D^h_{\w f(z)}\left(d_y \w F_c(w),v\right) d_y \D_z(u)\ d \mu^c_y(z)\\[1em]
\quad+\int_{\w Y} d_{\w F_c(y)} \D^h_{\w f(z)}(v)  \nabla d_y \D_z(w,u)\ d \mu_y^c(z)\\[1em]
\quad -c\int_{\w Y} d_{\w F_c(y)} \D^h_{\w f(z)}(v) d_y \D_z(u) d_y \D_z(w) \ d\mu^c_y(z)\\[1em]
\quad + c\, k'_y(u,v) \int_{\w Y} d_y \D_z(w) \ d\mu^c_y(z).
\end{array}
\end{equation}
Consider the second term in the right side of the previous equality. By condition  \eqref{ex(b)} the absolute values of the eigenvalues of the $\Hess \D_p$  are bounded by a positive constant $\lmb_0$, we have
\[
\int_{\w Y} d_{\w F_c\left( y\right) } \D^h_{\w f\left( z\right) }\left( v\right)\,   \nabla d_y \D_z\left( w,\, u\right) \, d \mu_y^c\left( z\right) \leq \lmb_0  \,\X\left( g\right)  \left\| v \right\| \left\| w \right\|\left\| u \right\|\int_{\w Y}\, d \mu_y^c.
\]
We can repeat a similar argument to any term of \eqref{epineq} and conclude that there exists constant $C>0$ such that \eqref{epineq0000} is verified. Analogously we can see that the integrands of the integrals in \eqref{epineq} and  \eqref{defnk'} are bounded, so that the previous derivations under the integral sign are  well defined.
\endproof

\begin{lem}\label{lemmaxxx}
For every $y \in \w Y \setminus \w Y_\ep$, with $\ep < \frac{1}{\left(2 B(n)\right)^4}$ and $c$ such that $\delta=\left(\frac{\X(g)\,c}{2\,n}\right)^{2n}-1< \frac{\ep}{(1+\ep)}$, we have
\begin{equation}\label{ineqlem9}
\begin{matrix}
\left| k'(u,v) - \frac{2}{c}\, \w g_{\hyp} \left(d_y \w F_c \left(u\right),v \right)\right|\leq\\[1ex]
\quad\leq B(n)\,\ep^{\frac{1}{4}}\, \frac{\left(B(n) \ep^{\frac{1}{4}} + \frac{2}{n}\right)^{\frac{1}{2}}\left(B(n) \ep^{\frac{1}{4}} + \frac{\X^2 (g)}{2\,n}\right)^{\frac{1}{2}}\left\| u \right\|\left\| v \right\|}{ \left(2 - B(n) \ep^{\frac{1}{4}}\right)} 
\end{matrix}
\end{equation}
for every $u \in T_y \w Y$, $v\in T_{\w F_c (y) } \C H^n$. 
\end{lem}
\proof 
By the definitions of $k'$, $K$ and equality \eqref{eq main1}, we have
\begin{equation*}
k'(u,v)= \frac{1}{c}\, \w g_{\hyp} \left( K \circ d_y \w F_c \left(u\right) , v \right)
\end{equation*}
hence
\[
\left| k'(u,v)- \frac{2}{c} \w g_{\hyp}\left(d_y \w F_c (u), v\right) \right|= \frac{1}{c}\, \left| \w g_{\hyp} \left(\left( K - 2I\right)d_y \w F_c (u), v\right) \right|
\]
\[
\leq \frac{1}{c} \left\| K-2I \right\| \left\| d_y \w F_c (u) \right\| \left\| v \right\|
\]
\[
\leq B(n)\,\ep^{\frac{1}{4}}\, \frac{\left(B(n) \ep^{\frac{1}{4}} + \frac{2}{n}\right)^{\frac{1}{2}}\left(B(n) \ep^{\frac{1}{4}} + \frac{\X^2 (g)}{2\,n}\right)^{\frac{1}{2}}\left\| u \right\|\left\| v \right\|}{ \left(2 - B(n) \ep^{\frac{1}{4}}\right)},
\]
where in the last inequality we used \eqref{epineq00} and \eqref{ineq k-2i}.
\endproof

\begin{lem}\label{main lemma}\label{epineq02}
If $\ep<\frac{1}{\left( {4\,B\left( n \right) }\right)^4 }$ and $c$ is such that $\delta=\left(\frac{\X(g)\,c}{2\,n}\right)^{2n}-1< \frac{\ep}{(1+\ep)}$, then for every $y \in Y$
\begin{equation}\label{epineq023}
\begin{split}
(1+\delta)^{2n}\xi(\ep)^{1-2n}\leq\frac{\left\| d_y \w F_c(u) \right\|}{\|u\|} \leq \xi(\ep),
\end{split}
\end{equation}
where
\begin{equation}\label{def xi}
\xi(\ep)=\frac{c\left(B(n) \ep^{\frac{1}{4}} + \frac{2}{n}\right)^{\frac{1}{2}}\left(B(n) \ep^{\frac{1}{4}} + \frac{\X^2 (g)}{2\,n}\right)^{\frac{1}{2}}}{ \left(2 - B(n) \ep^{\frac{1}{4}}\right)}
\end{equation}
\end{lem}
\proof
Suppose $\ep < \frac{1}{\left( 4\, B\left( n\right)  \right)^4}$. Let $H(y,r)=\Vol(Y,g)^{-1}  \int _{B\left( y,\, r\right)\subset Y}\,\nu_g$, due to the compactness of $Y$, it is a uniformly continuous map, so it is well defined the continuous function $h\left( r\right) =\min_y H\left( y,\, r\right) $.
Since $h\left( r \right) $ is strictly increasing, there exists an increasing function $\ep \f r(\ep)$ such that $h\left(r\left(\ep\right)\right)=2\,\sqrt \ep.$  By Lemma \ref{lemma ye} we see that
\begin{equation}\label{inclusion}
\begin{split}
B\left(y,r\left(\ep\right)\right) \not \subset Y_\ep
\end{split}
\end{equation}
for any $y$. Therefore, denoted $\w B\left(y,\, r\left(\ep\right)\right)=\pi^{-1} \left( B\left(y,\, r\left(\ep\right)\right)\right)$, we have  
\[
\w Y\setminus\w Y_\ep \cap \w B\left(y,\, r\left(\ep\right)\right)\neq \emptyset \qquad \forall \, y\in \w Y_{\ep}.
\]

By \eqref{inclusion}  for every $y''\in \w Y$ there exist $y \in \w Y \setminus \w Y_\ep$ such that the distance $d(y,y'')=r\leq r(\ep)$. Let $\gamma$ be  a minimizing geodesic with $\gamma(0)=y$ et $\gamma(r)=y''$. Set $\ep_0=\frac{1}{\left(4 B(n)\right)^4}$. We define $t_0 \in [0,r]$ the instant when $\gamma$ intersect $\w Y_{\ep_0}$ for the first time, if $\gamma$ does not intersect $\w Y_{\ep_0}$ we set $t_0=r$. So $\gamma \left([0,t_0] \right)\subset \w Y \setminus \w Y_{\ep_0}$. Define $y'=\gamma(t_0)$, let $u \in T_{y'}\w Y$ and $v \in T_{F_c(y')} \C H^n$, we define $U$ and $V$ the parallel field long $\gamma$ and $F(\gamma)$ such that $U(t_0)=u$ and $V(t_0)=v$. By Lemma \ref{lemHessK}
\begin{equation*}%\label{eq quasi}
\begin{split}
&\left|k'_{y'}(u,v)-k'_{y}\left( U(0), V(0) \right) \right| \leq  C\, d(y,y')\left( 1 + \sup_t \left\|d \w F \left( \dot \gamma(t) \right) \right\| \right)\left\|u\right\|\left\|v\right\|.
\end{split}
\end{equation*}
Therefore by \eqref{epineq00}, for any $y \in \w Y \setminus \w Y_{\ep} \text{ with } 0<\ep\leq\ep_0$ we have
\begin{equation*}%\label{eq quasi}
\begin{array}{C}
\left|k'_{y'}(u,v)-k'_{y}\left( U(0), V(0) \right) \right| \leq \\
C\, d(y,y')\left( 1 + \frac{c\left( \frac{1}{4} + \frac{2}{n}\right)^{\frac{1}{2}}\left(\frac{1}{4} + \frac{\X^2 (g)}{2\,n}\right)^{\frac{1}{2}}}{ \left(2 - \frac{1}{4}\right)}  \right)\left\|u\right\|\left\|v\right\|\leq\\[1em]
C\, r(\ep)\left( 1 + \frac{c}{7n} \left(n + 2\,\X^2 (g)\right)^\frac{1}{2} \left(n + 8\right)^\frac{1}{2} \right)\left\|u\right\|\left\|v\right\|,
\end{array}
\end{equation*}
hence, set $D\left( \ep \right) := C\, r(\ep)\left( 1 + \frac{c}{7n} \left(n + 2\,\X^2 (g)\right)^\frac{1}{2} \left(n + 8\right)^\frac{1}{2} \right)$, we get:
\begin{equation}\label{ineq k'}
\begin{split}
&k'_{y'}(u,v)\geq k'_{y}\left( U(0), V(0) \right)- D\left( \ep \right) \left\|u\right\|\left\|v\right\|.
\end{split}
\end{equation}
Since $V(0) \mapsto V(t_0)$ is an isometry between $T_{\w F_c(y)}\C H^n$ and $T_{\w F_c(y')}\C H^n$, and $\left| \Jac F_c \right|\left(\gamma(t)\right)\neq 0$, there exists $v\in T_{\w F_c(y')}\C H^n$, with $\left\| v \right\|=1$, such that $V(0)=v=\frac{d_{y'} \w F_c \left(U(0)\right)}{\left\| d_{y'} \w F_c \left(U(0)\right) \right\|}$. Let $K'_y:T_y \w Y \f T_{\w F_c(y)}\C H^n$ be the linear application defined by 
\[
\w g_{\hyp}\left( K'_y(u),w \right)=k'_y\left(u,w \right).
\]
By  \eqref{ineq k'}, we have
\begin{equation*}%\label{eq quasi}
\begin{array}{C}
\left\|K'_{y'}(u) \right\| \geq \w g _{\hyp} \left( K'_{y'}(u),v \right)=k'_{y'}\left( u,v \right)\geq k'_y\left( U(0), \frac{d \w F_c \left(U(0)\right)}{\left\| d \w F_c \left(U(0)\right) \right\|}\right)
- D\left( \ep \right) \left\|u\right\|.
\end{array}\end{equation*}
By \eqref{ineqlem9} and \eqref{def xi} we get
\begin{equation*}%\label{eq quasi}
\begin{split}
&\left\|K'_{y'}(u) \right\| \geq \frac{2}{c}\, \w g_{\hyp}\left(d \w F_c \left((U(0)\right), \frac{d \w F_c \left(U(0)\right)}{\left\| d \w F_c \left(U(0)\right) \right\|}\right)- B(n)\, \xi(\ep)\, \ep^{\frac{1}{4}}\left\| u \right\| - D\left( \ep \right) \left\|u\right\|.
\end{split}\end{equation*}
By \eqref{epineq01} we obtain
\begin{equation*}%\label{eq quasi}
\begin{split}
&\left\|K'_{y'}(u) \right\| \geq \left(\frac{2\left( \left( 1+\delta\right)\left( 1- \sqrt \ep \right)   \right)^{{2\,n}}}{c\, \xi^{2n-1}(\ep)}- B(n)\, \xi(\ep)\, \ep^{\frac{1}{4}} - D\left( \ep \right) \right)\left\| u \right\|,
\end{split}\end{equation*}
so
\begin{equation}\label{detK'}
\begin{split}
&\left|\det \left(K'_{y'}\right)\right| \geq \left(\frac{2\left( \left( 1+\delta\right)\left( 1- \sqrt \ep \right)   \right)^{{2\,n}}}{c\, \xi^{2n-1}(\ep)}- B(n)\, \xi(\ep)\, \ep^{\frac{1}{4}}- D\left( \ep \right) \right) ^{2n},
\end{split}\end{equation}
By the definitions of $k'$ and $K$ and equality \eqref{eq main1}, we have
\begin{equation*}
\w g_{\hyp} \left( K \circ d_y \w F_c \left(u\right) , v \right)=c\, k'_y(u,v)=c\,\w g_{\hyp}\left( K'_y(u),v \right)
\end{equation*}
thus
\begin{equation*}
\det\left(K_{y'} \right) \Jac \w F_c\left( y'\right)=c^{2n}\det\left( K'_{y'} \right),
\end{equation*}
by \eqref{hessian compl hyperb} we see that $\Tr K_{y'}= 4n$. So we get
\begin{equation*}%\label{eq quasi}
\begin{split}
&\left| \Jac \w F_c \right|\left( y'\right)= \left( \frac{\Tr K_{y'}}{4n} \right)^{2n}\left| \Jac \w F_c \right|\left( y'\right)\\
&\quad \geq\left(\frac{1}{2}\right)^{2n}\left|\det \left( K_{y'} \right) \right|\left| \Jac \w F_c \right|\left( y'\right) \\
&\quad=\left(\frac{c}{2}\right)^{2n}\left|\det \left( K'_{y'} \right) \right|\left( y'\right)
\end{split}\end{equation*}
therefore by \eqref{detK'}
\begin{equation}\label{boundF}
\begin{split}
&\left| \Jac \w F_c \right|\left( y'\right)\geq c^{2n} \left(\frac{\left( \left( 1+\delta\right)\left( 1- \sqrt \ep \right)   \right)^{{2\,n}}}{c\, \xi^{2n-1}(\ep)}- \frac{1}{2} B(n)\, \xi(\ep)\, \ep^{\frac{1}{4}}- \frac{1}{2}\,D\left( \ep \right) \right)^{2n}  .
\end{split}\end{equation}
If $\gamma$ intersect $\ov{\w Y_{\ep_0}}$, by \eqref{boundF}, we have
\begin{equation*}%\label{boundF}
\begin{split}
&(1-\sqrt {\ep_0}) (1+\delta)\geq c^{2n}  \left(\frac{\left( \left( 1+\delta\right)\left( 1- \sqrt \ep \right)   \right)^{{2\,n}}}{c\, \xi^{2n-1}(\ep)}- \frac{1}{2} B(n)\, \xi(\ep)\, \ep^{\frac{1}{4}}- \frac{1}{2}\,D\left( \ep \right) \right)^{2n}  .
\end{split}\end{equation*}
Since the previous inequality hold for $0<\ep \leq \ep_0$ and $\delta< \frac{\ep}{(1+\ep)}$, we get a contradiction as $\ep$ approach to zero, indeed the first member goes to $(1-\sqrt {\ep_0})$, on the contrary the second member  goes to $1$. We conclude that $\w Y_{\ep_0}=\emptyset$. Therefore, passing $\w F_c$ to its quotient $F_c$, equations \eqref{epineq00} and \eqref{epineq01} imply \eqref{epineq023}.
\endproof

\noindent{\bf Proof of Theorem \ref{thm almost}:} 
Set 
$
A'(\ep)=\left( 1+\delta\right) ^{2\,n}\xi\left( \ep\right) ^{1-2\,n}
$
and
$
A''\left( \ep\right) =\xi\left( \ep\right) 
$
in Lemma \ref{epineq02},
where $\xi(\ep)$ is given by  \eqref{def xi} (notice that $\xi(\ep) \geq 1$). \qed

\smallskip

{{The proof of the first part of Theorem \ref{thm dentropy2} is complete.}}

\medskip

{\noindent{\bf Conclusion of the proof of Theorem \ref{thm dentropy2}, the $\ep=0$ case.} We want to prove that when $\ep = 0$}, then $F_c$ is a holomorphic or anti-holomorphic local isometry.   Suppose that $g$ is normalized in order to have
\[
\Vol\left( Y,\, g\right)  =\deg\left( f \right)  \Vol\left( X,\, g_{\hyp}\right),
\]
we want to prove that there exists a riemannian covering $F:\left(Y,g\right)\f \left(X,g_{\hyp}\right)$.\\
\noindent Take a sequence $\left\{ F_{c_n}\right\}$ such that $0<\X(g)\,c_n - c_0< c_0\left( \sqrt[2n]{\frac{1+2\frac{1}{n}}{1+\frac{1}{n}}} -1\right) $. For $n$ sufficiently large, say $n>n_0$, the sequence $\left\{ F_{c_n} \right\}$ consists of $\mathcal C^1$ covering maps. Being $X$ and $Y$ compact the $F_{c_n}$ are  equibounded. By inequalities \eqref{eq quasi} we get
\begin{equation*}%\label{eq quasi}
\begin{split}
&\left\|F_{c_n}(y_0) - F_{c_n}(y) \right\| \leq  \left\|d_{y_0} F_{c_n}\right\|\left\|y_0-y\right\|\leq A''(\ep)\left\|y_0-y\right\|\\
&\quad\leq A''\left(\frac{1}{n}\right)\left\|y_0-y\right\|,
\end{split}\end{equation*}
therefore the maps $F_{c_n}$ are equicontinuous. By the Ascoli-Arzel\`a theorem there exist a subsequence $c_n \f c_0$, such that %(with abuse of notation)
$F_{c_n}$ uniformly converge to a continuous function $F$ with $\deg\left( F \right) = \deg\left( F_{c_n}\right) = \deg\left( f \right)$. 
Let $\gamma:[0,1]\f Y$ a piecewise regular curve such that $\gamma(0)=y_1$ and $\gamma(1)=y_2$  then
\[
\int_0^1 A'\left( \frac{1}{n}\right)  \left\| \dot \gamma (t)\right\| \, dt \leq \int_0^1 \left\|d_y F_{c_n}\left( \dot \gamma \left( t\right) \right) \right\| \, dt \leq \int_0^1 A''\left( \frac{1}{n}\right) \left\|\dot \gamma \left( t\right) \right\| \, dt
\]
hence, denoted respectively $d\left( \cdot , \cdot \right)$  and $d^h\left( \cdot , \cdot \right)$ the geodesic distance on $\left( Y,\, g \right) $ and $\left( X ,\, g_{\hyp} \right)  $  we have
\[d^h \left(F\left(y_1\right),\, F \left(y_2\right)\right)\leq d \left(y_1,\, y_2\right).
\]
By \cite[Proposition C.1]{bcg1} the map $F$ is a riemannian covering.
Arguing as in the last part of proof of \cite[Theorem 1.1]{M06} we deduce that $F$ is holomorphic or anti-holomorphic. The proof of Theorem \ref{thm dentropy2} is complete.

\vskip 0.3cm

\noindent {\bf Proof of Theorem \ref{thm proper}}: we need to verify that   conditions \eqref{ex(a)} and \eqref{ex(b)} above are satisfied and then apply Theorem \ref{thm dentropy2}.

\vskip 0.3cm

\noindent {\bf Condition \eqref{ex(a)} is satisfied.}  
 Let $\varphi: (Y, \, g) \f (N, \, g^{N})$ be the strongly proper \K\ immersion of $Y$ in an locally classical symmetric space of noncompact type $N$ and let $ \w \varphi: (\w Y, \, \w g) \f (\Omega, \,  g^\Omega)$ be its lift to the \K\ universal covers. By \eqref{hered glob diast} we see that $(\w Y, \, \w g)$ has the diastasis globally defined.
As $\lim_{t \f +\infty}\frac{\log\cosh t} { t} = 1$, fixed $\mu >0$ and $q \in \w Y$, there exists a compact set $K\subset \w Y$, two constant $L_{1}, \, L_{2} \in \R $ and $\varepsilon > 0$ such that $\forall \, p \in \w Y \setminus  K $, 
\begin{equation}\label{int matrix 1}
\begin{matrix}
\rho \left(q,\, p\right) - L_2 <  L_1 \, e^{\,\mu\, \rho^{\Omega} \left(\w\varphi\left(q\right),\, \w\varphi\left(p\right)\right)} 
< L_1 \, e^{\,\mu\, \left( \log \cosh  \rho^{\Omega} \left(\w\varphi\left(q\right),\, \w\varphi\left(p\right)\right) + \varepsilon \right)} \\[1em]
< L _1\, e^{\,\frac{\mu}{2}\, \left( \D^{\Omega} \left(\w\varphi\left(q\right),\, \w\varphi\left(p\right)\right)\, +\, 2\,\varepsilon \right)}=
L _1 \, e^{\,\frac{\mu}{2}\, \left( \D \left( q,\, p \right)\, +\, 2\,\varepsilon \right)}, 
\end{matrix}
\end{equation}
where in the first inequality we use that $\varphi$ is strongly proper (notice that this is the unique point of the proof where this hypothesis is used), in the third one we used \eqref{diast poly e distantance}, while in the last equality we applied \eqref{hered glob diast}. 
On the other hand, if we choose $\frac{\mu}{2} >0$ small enough so that $(c-\frac{\mu}{2}) > \frac{\Ent\left( Y,\, g\right)}{\X \left( g \right)}$ we obtain
\begin{equation}\label{int matrix 2}
\int_{\w Y\setminus K } e^{-(c-\frac{\mu}{2})\, \D_q\left( p\right) }\,\nu_{g}\left( p \right) < \infty.
\end{equation}
Putting together \eqref{int matrix 1} and \eqref{int matrix 2} we see that
$
\int_{\w Y}  \rho \left( q, \, p \right) e^{-c\, \D_q\left( p\right) }\,\nu_{g}\left( p \right)
$
is convergent, so \eqref{ex(a)} is verified. 

\vskip 0.3cm

\noindent {\bf Condition \eqref{ex(b)} is satisfied.}
Being $Y$ compact, the second fundamental form of $\w \varphi$ is bounded.  Hence the conclusion follow by combining \eqref{HessDiast},  \eqref{grad class thm} and \eqref{hessian bounded} setting $\psi = \w \varphi$.
\hfill \ensuremath{\Box}
\endproof

%\vskip 0.3cm
%\noindent {\bf Observation.} Notice that with the same argument we can actually prove that if the universal \K\ covering of a compact \K\ manifold $\left(Y,\, g\right)$  admit a \K\ immersion $\w\varphi : \left(\w Y,\, \w  g\right) \f \left(\Omega,  g^{\Omega}  \right)$ in a classical HSSNCT such that
%\begin{enumerate}[(i)]
%\item Fixed $\mu > 0 $ and $q\in \w M_1$, there exist two constants $L_1>0$ and $L_2>0$, such that
%$$\rho \left(q,\, p\right) < L_1 \, e^{\,\mu\, \rho^\Omega \left(\varphi\left(q\right),\, \varphi\left(p\right)\right)} + L_2, \qquad \forall \, p \in \w Y$$
%where $\rho$ and $\rho^{\Omega}$ are the geodesics distances on $\w Y$ and $\Omega$ respectively;
%\item the second fundamental form of $\w\varphi$ is bounded, $\sup_{y \in \w Y} \left\| {\operatorname{II}_y}\right\|<\infty$;
%\end{enumerate}
%then \eqref{ex(a)} and \eqref{ex(b)}  are verified.
%

\end{document}